\input amstex
\documentstyle{amsppt}
\input bull-ppt
\keyedby{bull297e/pah}

\topmatter
\cvol{27}
\cvolyear{1992}
\cmonth{July}
\cyear{1992}
\cvolno{1}
\cpgs{148-153}


\title Rational Function Certification\\
Of Multisum/integral/``$q$'' Identities
\endtitle
\shorttitle{Multivariate hypergeometric identities}
\author Herbert S. Wilf and Doron Zeilberger\endauthor
\shortauthor{H. S. Wilf and Doron Zeilberger}
\address University of Pennsylvania, Philadelphia, 
Pennsylvania 19104-6395\endaddress
\ml wilf\@central.cis.upenn.edu\endml
\address Temple University, Philadelphia, Pennsylvania 
19122-2585\endaddress 
\ml zeilberg\@euclid.math.temple.edu\endml
\date December 2, 1991 and, in revised form,  January 28, 
1992\enddate
\thanks 
The first author was supported in part by the U. S. Office 
of Naval Research\endthanks 
\thanks The second author
was supported in part by the U. S. National Science 
Foundation\endthanks
\abstract The method of rational function certification 
for proving terminating
hypergeometric identities is extended from  single sums or 
integrals to
 multi-integral/sums and ``$q$'' integral/sums.\endabstract

\endtopmatter
\document

\heading 1. Introduction\endheading

{\it Special functions} have been defined by Richard Askey 
[As2] as
``functions that occur often enough to merit a name,'' 
while Tur\' an
[As1, p.\ 47] defined them as ``useful functions.'' The 
impact of
these special functions on classical mathematics and 
physics can be
gauged by the stature of those whose names they bear: 
Bessel, Gauss,
Hermite, Jacobi, Legendre, Tschebycheff, to name a few. It 
turns out that most
special functions are of {\it hypergeometric type}, which 
is to say
that they can be written as a sum in which the summand is a 
hypergeometric term. Also of great interest are the 
so-called $q$-analogs
of special functions and hypergeometric series, called 
$q$-series.
These have many applications to number theory, 
combinatorics, physics,
group theory [An], and other areas of science and 
mathematics.

There are countless identities relating special functions 
(e.g.,
[PBM, R, An, As1]). In addition to their intrinsic interest,
some of them imply important properties of these special 
functions,
which in turn sometimes imply deep theorems elsewhere in 
mathematics
(e.g., [deB, Ap]). Just as important for mathematics are the
extremely successful attempts to instill meaning and 
insight, both
representation-theoretic (e.g., [Mi]) and combinatorial 
(e.g., [Fo]),
into these identities.

The general theory of multivariate hypergeometric 
functions is
currently a very active field, rooted in multivariate 
statistics and
the physics of angular momentum. A very novel and fruitful 
approach is 
currently being pursued by Gelfand, Kapranov, Zelevinsky 
(e.g., [GKZ]) and
their collaborators.

We now know [Z1, WZ1] that terminating identities 
involving sums and
integrals of products of special functions of 
hypergeometric type can
be proved by a finite algorithm, viz., find recurrence or 
differential
equations that are satisfied by the left and the right 
sides of the
claimed identity (they always exist), and then compute 
enough initial
values of the two sides to assure that the two recurrence or
differential equations have the same solution. We know 
further that
for {\it one-variable} hypergeometric sums [WZ1, WZ2, Z2, 
Z3] and
single ``hyperexponential'' integrals [AZ] there are {\it 
efficient}
algorithms for doing this, so computers really can do the 
job.

Here we announce fast, efficient algorithms for such 
identities that
involve {\it multiple} sums and integrals of products of 
special functions
of hypergeometric type. We also announce the algorithmic 
provability
of single- and multivariate $q$-identities and, 
furthermore, that the
algorithms are fast and efficient. The proofs that are 
generated by
these programs are extremely short, and human beings can 
verify them
easily.

Further we show how to extend to ``explicitly-evaluable'' 
{\it
multiple} sums and integrals the notions of {\it rational 
function
certification}, {\it WZ-pair}, and {\it companion 
identity} that were
introduced in [WZ1] and [WZ2] for single sums.

The present theory applies directly only to terminating 
identities, by
which we mean that for any specific numerical assignment 
of the
auxiliary parameters the integral-sum is trivially 
evaluable. However,
very often {\it nonterminating} identities are limiting 
cases or
``analytic continuations'' of terminating ones, so our 
results also
bear on them.

 Our method proves (or refutes) any such given conjectured
 identity, but humans are still needed to conjecture 
interesting ones.

Full details and many examples will appear elsewhere 
[WZ3]. Our
Maple programs are available from $<${\tt 
zeilberg\@euclid.math.temple.edu}$>$.

\heading 2. Hypergeometric multivariate 
identities\endheading

For a function $F(\bold k,\bold y)$ of $r$ discrete 
variables ${\bold
k}$ and $s$ continuous variables ${\bold y}$ let $K_i$ be 
the operator
defined by $K_iF(\bold k,\bold y)=F(k_1,\dots ,k_i+1,\dots 
,k_r,\bold
y)$, and let $D_j=\partial /\partial y_j$, for $i=1,\dots 
,r$ and
$j=1,\dots ,s$.

\proclaim{Definition} A function $F(\bold k,\bold y)$ is a 
hypergeometric term
 if $K_iF/F$
$(i=1,\dots ,r)$ and $D_jF/F$ $(j=1,\dots ,s)$ are all
 rational functions of $(\bold k,\bold y)$.\endproclaim

A sequence of special functions of hypergeometric type (in 
one
variable) is a sequence that is of the form 
$P_n(x):=\sum_kF(n,k)x^k$
where $F$ is a hypergeometric term.

A typical identity in the theory of special functions 
involves
multiple integral-sums of products of polynomials of 
hypergeometric
type. After a full expansion, such an identity is of the 
form {\rm ``}left
side{\rm ''}={\rm ``}right side{\rm ''} 
where both sides are of the form $\sum_{\bold
k}\int_{\bold y}F(\bold k,\bold n,\bold x,\bold y)d{\bold 
y}.$ Here
$F$ is a hypergeometric term in the discrete 
multivariables $(\bold
k,\bold n)$ and the continuous multivariables $(\bold 
x,\bold y)$. It
might happen that one of the sides, say the right side, 
has no
$\sum$'s or $\int$'s in it, i.e., it is already 
hypergeometric, in
which case we speak of {\it explicit evaluation}.

 Finally, we also require that our integrand/summands be 
holonomic.
This is certainly the case for what we call 
proper-hypergeometric
terms, which for purely discrete summands, look like 
$$F(n,\bold
k)={{\prod (a_in+{\bold b}_i\cdot {\bold k}+c_i)!}\over 
{\prod
(u_in+{\bold v}_i\cdot {\bold k}+w_i)!}}P(n,\bold k 
){\boldxi}
^{\bold k},\eqno(\text{P-H})$$ where the $a$'s and $u$'s 
are specific
 integers, ${\bold b}$ and ${\bold v}$ are vectors of 
specific integer
entries, the $c$'s and the $w$'s are complex numbers that 
may depend
upon additional parameters, $P$ is a polynomial in ${\bold 
k}$, and
${\boldxi}$ is a vector of parameters.

This allows us to use the holonomic theory of [Z1].
However for pure multisums, we give proofs of our results,
with effective bounds, that are self-contained and
independent of the theory of holonomic functions.

\heading 3. The fundamental theorem\endheading

Let $\Delta_i=K_i-1$ be the forward difference operator in 
$k_i$, and
let $N$ be the forward shift in $n$: $Nf(n)=f(n+1)$.

\proclaim{Theorem 1} Let $F(n,\bold k,\bold y)$ \RM(resp.\ 
$F(x,\bold k,\bold y))$
 be a proper-hypergeometric term, or more generally, a 
holonomic
hypergeometric term, in $(\bold k,\bold y)$ and $n$ 
\RM(resp.\ $x)$, where
$n,\bold k$ are discrete and $x,\bold y$ are continuous 
variables.
Then there exist a linear ordinary recurrence \RM(resp.\ 
differential\,\RM)
operator $P(N,n)$ \RM(resp.\ $P(D_x,x))$ with polynomial 
coefficients and
rational functions $R_1,\dots ,
R_r, S_1,\dots ,S_s$ such that
$$P(N,n)F\quad (\text{resp.
}P(D_x,x)F)=\sum_{i=1}^r\Delta_{k_i}(R_iF)+
\sum_{j=1}^sD_{y_j}(S_jF).\tag1$$
\endproclaim

For proper-hypergeometric $F$ we give explicit a priori 
bounds for the
order of the operator $P$. By summation-integration of (1) 
over $\bold
k,\bold y$ we obtain the following

\proclaim{Corollary A} If $F(n,\bold k,\bold y)$ 
\RM(resp.\ $F(x,\bold k,
\bold y))$
 is as above and of compact support in $(\bold k,\bold y)$ 
for every
fixed $n$ \RM(resp. $x)$ then $$f(n)\ (\text{resp. 
}f(x)):=\sum_{\bold
k}\int_{\bold y}Fd\bold y\tag 2$$ satisfies a linear 
recurrence
\RM(resp.\ differential\,\RM) equation with polynomial 
coefficients
$$P(N,n)f(n)\equiv 0\quad (\text{resp. } P(D_x,x)f(x)\equiv
0).\tag 3$$\endproclaim

The denominators of the rational functions $\bold R,\bold 
S$ can be
predicted for any given $F$, and an upper bound for the 
order of the
operator $P$ can be given in advance. Hence by assuming 
the operator
and the rational functions in the most general form with 
those
denominators and with that order, the determination of the 
unknown
operator and rational functions quickly reduces to solving 
a system of
linear equations with {\it symbolic} coefficients. By [W] 
this reduces
to solving such a system with {\it numerical} 
coefficients, and that
in turn reduces to solving linear systems with {\it integer}
coefficients, a problem for which fast parallelizable 
algorithms exist
[CC].

Using the terminology of [WZ1, WZ2] one can say that the 
rational
functions $(\bold R,\bold S)$ {\it certify} the recurrence 
(resp.\ the
differential equation) (3).

\proclaim{Corollary B} Any identity of the form 
{\rm ``}left side{\rm ''}={\rm ``}right side,{\rm ''} 
where both sides have the form {\rm (2)} and $F$ is a 
proper-hypergeometric
or holonomic-hypergeometric term, has a two-line 
elementary proof,
constructible by a computer and verifiable by a human or a
computer.\endproclaim

\heading 4. Sketch of the proof of the theorem\endheading

From the general theory of [Z1] we can find a linear partial
recurrence-differential operator $T(n,N,\bold K,\bold 
D_{\bold y})$,
independent of $\bold k,\bold y$, that annihilates $F$. For
proper-hypergeometric functions we have a completely 
elementary proof
that also gives explicit upper bounds for the orders in 
$N$ and $\bold
K$. This elementary proof was obtained by extending to
multisums-integrals Sister Celine Fasenmyer's technique 
[Fa, R] as
systematized by Verbaeten [V].

Once the operator $T$ has been found, we write 
$$T(n,N,\bold K,\bold
D_{\bold y})=P(N,n)+\sum_{i=1}^r(K_i-1)T_i(N,\bold K,\bold
D)+\sum_{j=1}^sD_j\widehat{T}_j(N,\bold K,\bold D).\tag 
4$$ It is
easy to see that this is always possible. Next set
$$\aligned
G_i(n,\bold k,\bold y)&:=-T_iF(n,\bold k,\bold y)\qquad
(i=1,\dots ,r),\\
\widehat{G}_j(n,\bold k,\bold y)&:=-\widehat{T}_jF(n,\bold 
k,\bold y)\qquad 
(j=1,\dots ,s).\endaligned\tag 5$$

Since $F$ is hypergeometric, the $G_i$ and the 
$\widehat{G}_j$ are
rational multiples of $F$: $G_i=R_iF$, 
$\widehat{G}_j=S_jF$. Now apply
(4) to $F(n,\bold k,\bold y)$, remembering that $TF=0$, to 
get (1). \qed 

\heading 5. Discrete and continuous $q$-analogues\endheading

The above extends to multivariate $q$-hypergeometric 
identities. Let
$Q_j$ be the operator that acts on $F$ by replacing $y_j$ 
by $qy_j$
wherever it appears. Then we say that a function $F(\bold 
k,\bold y)$
is a $q$-hypergeometric term if for each $i=1,\dots ,r$ 
and $j=1,\dots ,s$,
it is true that $K_iF/F$ and $Q_jF/F$ are rational 
functions of
$(q,q^{k_1},\dots ,q^{k_r},y_1,\dots ,y_s)$. There is also 
a natural
definition of $q$-{\it proper-hypergeometric}, which is 
given in
[WZ3].

The fundamental theorem still holds, where the integration 
is replaced
either by Jackson's $q$-integration [An] or by an ordinary 
contour
integral, or, in the case of a formal Laurent series, by the
action of taking ``constant term of.'' Macdonald's 
$q$-constant term
conjectures ([Ma], see [Gu, GG] for a recent update) for 
every
{\it specific} root system, fall under the present heading.

\heading 6. Explicit closed-form identities:\\
WZ-tuples and companion identities\endheading

In the case where the identity ``left side''=``right 
side'' is such
that the right side does not contain any `$\sum$' or 
`$\int$' signs,
i.e., is of closed form, one has an {\it explicit 
identity}. If the
right side is nonzero one can divide through by it to get 
an identity
of the form $$\sum_{\bold k}\int_{\bold y}F(n,\bold 
k,\bold y)\,d\bold
y=1.$$ Since the summand satisfies (1), the left side, 
call it $L(n)$,
satisfies some linear recurrence $P(N,n)L(n)=0$, by 
Corollary A. Often
the operator $P$ turns out to be the minimal order 
recurrence that is
satisfied by the sequence that is identically 1, viz. 
$(N-1)L(n)=0$.
If that happens then if we let $G_i:=-R_iF$ and 
$H_j:=-S_jF$ we find
that (1) becomes
$$
\Delta_nF+\sum_{i=1}^r\Delta_iG_i+\sum_{j=1}^sD_{y_j}H_j=0.
\tag6
$$ 
We call $(F,\bold G,\bold H)$ a {\it WZ-tuple}. It 
generalizes the
idea of {\it WZ-pair} developed in [WZ1, WZ2]. Recall that 
a WZ pair
gave, as a bonus, a new identity, the {\it companion 
identity}. Here,
if we sum-integrate (6) w.r.t. all of the variables except 
{\it one},
we get a new identity for each choice of that one 
variable, for a
total of $r+s$ new companion identities altogether!

\heading 7. Example: The Hille-Hardy Bilinear Formula\\
 for Laguerre Polynomials
\endheading
 
 As an example we will now show the computer proof of the 
Hille-Hardy
formula [R, Theorem 69, p. 212], namely, 
$$
\align
\noalign{${{n!}\over {( \alpha +1)_n }}L_n^{(
\alpha )} (x)L_n^{( \alpha )}(y)=
{1\over {2\pi i}}\int_{|u| = \epsilon } u^{-n-1} 
(1-u)^{-\alpha -1}  
\exp\left\{-{{(x+y)u}\over {(1-u)}}\right\}$}\\
\shoveright{\times  
\( \sum_m {i \over {m!( \alpha +1)_m}}\left( {{xyu} \over 
{(1-u)^2}}
\right)^m\)\, du.}
\endalign
$$  
Many other examples appear in \cite{WZ3}.  

To this end, it is enough to prove that the right side is 
annihilated
by the well-known second-order differential operator 
annihilating the
Laguerre polynomials, both w.r.t. $x$ and $y$. Of course, 
by symmetry,
it suffices to do it only for $x$, but the computer does 
not mind doing
it for both $x$ and $y$. We still need to prove that the 
initial
conditions match, but they are just the usual defining 
generating
function for the Laguerre polynomials.  The computer 
output was as
follows.
  
\proclaim{Theorem} Let 
$$F(u,m,x):={{(1-u)^{-\alpha -1}\exp{( -(x+y)u/(1-u)) 
(xyu/( 1-u)^2 )^m }} 
\over { u^{n+1}m!\Gamma ( \alpha +1+m)}},$$ 
and let $a(x)$ be its contour integral w.r.t. $u$   and 
sum w.r.t.   $m$.   
Let $D_x$ be differentiation w.r.t. $x$.  The function 
$a(x)$ satisfies  
the differential equation $$(n+(\alpha +1-x) D_x + x D_x^2 
) a(x) =0.$$  
\endproclaim
  
\noindent{Proof.} It is routinely verifiable that 
$$\aligned
&(n + ( \alpha + 1 - x) D_x + x D_x^2 ) F(u, m, x) \\
&\qquad= D_u( - u F(u, m,
x))
+\Delta_m(-(m ( \alpha + m)/x) F(u, m, x))
\endaligned$$ 
and the result follows
by integrating w.r.t. $u$ and summing w.r.t.  $m$. \qed

\rem{Remark} The phrase ``routinely verifiable'' in the 
above means
 that after carrying out the indicated differentiation and
differencing, and after dividing through by $F$ and clearing
denominators, what will remain will be a trivially 
verifiable {\it
polynomial} identity.
\endrem

\heading Acknowledgment\endheading
Many thanks are due to Dr. James C. T. Pool,
 head of Drexel University's Department of Mathematics and 
Computer
Science, for his kind and generous permission to use his 
department's
computer facilities.

\enddocument